\documentclass[1p]{elsarticle}
\usepackage{amssymb}
\usepackage{amsfonts}
\usepackage{amsmath}
\usepackage{pgf,tikz}
\usetikzlibrary{arrows}
\usetikzlibrary{patterns}
\usepackage{caption}

\newtheorem{theorem}{Theorem}
\newtheorem{conjecture}[theorem]{Conjecture}
\newtheorem{corollary}[theorem]{Corollary}

\newtheorem{proposition}[theorem]{Proposition}
\newtheorem{claim}{Claim}

\newcommand{\BIN}{\ensuremath{\operatorname{BIN}}}

\newproof{pf}{Proof}

\begin{document}
\title{On inclusion chromatic index of a graph}

\author[AGH]{Jakub Kwa\'sny\fnref{MNiSW}}
\ead{	}

\author[AGH]{Jakub Przyby{\l}o\corref{cor1}\fnref{MNiSW}} 
\ead{jakubprz@agh.edu.pl, phone: 048-12-617-46-38,  fax: 048-12-617-31-65}

\cortext[cor1]{Corresponding author}
\fntext[MNiSW]{This work was partially supported by the Faculty of Applied Mathematics AGH UST statutory tasks within subsidy of Ministry of Science and Higher Education.}

\address[AGH]{AGH University of Science and Technology, al. A. Mickiewicza 30, 30-059 Krakow, Poland}

\begin{abstract}
Let $\chi'_\subset(G)$ be the least number of colours necessary to properly colour the edges of a graph $G$ with minimum degree $\delta\geq 2$ so that the set of colours incident with any vertex is not contained in a set of colours incident to any its neighbour. We provide an infinite family of examples of graphs $G$ with $\chi'_\subset(G)\geq (1+\frac{1}{\delta-1})\Delta$, where $\Delta$ is the maximum degree of $G$, and we conjecture that $\chi'_\subset(G)\leq \lceil(1+\frac{1}{\delta-1})\Delta\rceil$ for every connected graph with $\delta\geq 2$ which is not isomorphic to $C_5$. The equality here is attained e.g. for the family of complete bipartite graphs. Using a probabilistic argument we support this conjecture by proving that for any fixed $\delta\ge2$, $\chi'_\subset(G) \le (1+\frac{4}{\delta})\Delta (1+o(1))$ (for $\Delta\to\infty$), what implies that $\chi'_\subset(G) \le (1+\frac{4}{\delta-1})\Delta$ for $\Delta$ large enough.
\end{abstract}

%\subjclass{05C15}
%\subjclass[2010]{05C15}

\maketitle

\section{Introduction}

Let $G=(V,E)$ be a (simple and finite) graph and $c\colon E\rightarrow C$ an edge colouring of $G$. For each vertex $v\in V$ we define a \emph{palette} of $v$ as
$$S_c(v) = \{c(uv): u\in N(v)\}.$$
This shall be referred to simply as $S(v)$ if it is clear which colouring is to be considered. We call $c$ an \emph{inclusion-free} colouring if it is proper, i.e. attributes distinct colours to adjacent edges, and for each $uv\in E$, $S_c(u)\not\subseteq S_c(v)$. Note that this condition may be restricted to the cases when $d_G(u)\le d_G(v)$ exclusively. Such a colouring exists if and only if $G$ has no vertices of degree one. The least number of colours admitting an inclusion-free colouring of $G$ is called its \emph{inclusion chromatic index} and denoted by $\chi'_\subset(G)$. This problem was proposed by~Simonyi~\cite{SimonyiMalta}, inspired by the concept discussed in the following paragraph (and in particular its list analogue, see~\cite{KwasnyPrzybylo_list_Hatami}).

The inclusion-free colouring is in fact a strengthening of so-called \emph{adjacent vertex distinguishing colouring} which requires the palettes of adjacent vertices to be merely distinct -- the minimizing necessary number of colours graph invariant is in this case denoted by $\chi'_a(G)$. In particular, the two problems are equivalent when $G$ is a regular graph. Moreover, by the properness of edge colourings investigated, the both corresponding graph invariants equal at least $\Delta$ or $\Delta+1$ (dependent on the class of a graph), while it was boldly conjectured by Zhang et al.~\cite{Zhang} that $\chi'_a(G)\leq \Delta+2$ for any connected graph $G$ of order at least 3 which is not the cycle $C_5$. This problem was widely studied \cite{Akbari,BalGLS,BonamyEtAl,BuLihWang,Hatami,HocqMont,Hornak_planar,LiZhangChenSun,WangWang,Zhang}. In particular, Balister et al. proved the conjecture to hold for bipartite and for subcubic graphs, and showed that in general $\chi'_a(G)\leq \Delta+O(\log\chi(G))$. Later Zhang's Conjecture was proved up to an additive constant by Hatami~\cite{Hatami}, who applied probabilistic approach to show that $\chi'_a(G)\leq \Delta+300$ for any graph $G$ with no isolated edges and with maximum degree $\Delta > 10^{20}$. Recently, Joret and Lochet \cite{Joret} used the so-called \emph{entropy compression method} and made yet another step forward towards the conjecture, improving Hatami's bound to $\chi'_a(G)\leq\Delta+19$ (for large enough $\Delta$). 

In this paper we first argument that contrary to our initial supposition (based on resemblance in the case of regular or almost regular graphs), the two problems are not much alike. Namely, in the next section we discuss the family of examples, the complete bipartite graphs, proving that no upper bound of the form $\Delta+const.$ can be expected in the case of $\chi'_{\subset}(G)$ in general, as we may have $\chi'_{\subset}(G)\geq (1+\frac{1}{\delta-1})\Delta$, where $\delta$ is the minimum degree of $G$. We then pose in Section~\ref{SectionConjecture} a conjecture that such quantity rounded up to an integer is essentially also an upper bound for $\chi'_{\subset}(G)$, and support this conjecture by proving that for every fixed $\delta\geq 2$ it holds up to an additive constant with $\frac{1}{\delta-1}$ replaced by $\frac{4}{\delta-1}$, see Section~\ref{ProbabilisticProofInclusion} for details of a probabilistic proof of this fact. This is also preceded by a useful deterministic argument that $\chi'_{\subset}(G)< 3\Delta$.

\section{Complete Bipartite Graphs}

The first difference between the two problems is the already mentioned fact that investigating inclusion-free colourings makes sense only for graph without degree one vertices (while we need only exclude isolated edges in the case of $\chi'_a(G)$), and thus from now on we assume that every graph considered has the minimum degree $\delta\geq 2$. However, even in such a case the influence of $\delta$ on $\chi'_{\subset}(G)$ can still be significant,
especially when $\delta$ is small compared to $\Delta$, as exemplified by the following result for complete bipartite graphs.

\begin{proposition}\label{ThBipartiteComplIncl}
Let $G=K_{\Delta,\delta}$, $2\leq \delta\leq \Delta$. Then
$$\chi'_{\subset}(G) = \left\lceil\left(1+\frac{1}{\delta-1}\right)\Delta\right\rceil.$$
\end{proposition}

\begin{pf} 
Let $G=(X\cup Y, E)$ where $|X|=\delta$, $|Y|=\Delta$, and consider its inclusion-free colouring $c$ with elements of $C$, $|C|=\Delta+k$. Let $x\in X$ and denote by $G'$ the subgraph induced by the edges coloured with colours within $C':=C\smallsetminus S(x)$ in $G$ (where $|C'|=k$). As $c$ is inclusion free, every vertex from $Y$ must be incident with an edge coloured by a colour in $C'$, and hence $G'=(X'\cup Y,E')$ is a properly coloured bipartite graph with at least $\Delta$ edges, $|X'|\leq \delta-1$ (as $x\notin X'$) and $d_{G'}(x')\leq k$ for every $x'\in X'$. Therefore, $k\ge \lceil\frac{\Delta}{\delta-1}\rceil$, what yields the desired lower bound for $\chi'_{\subset}(G)$.

For integers $n$ and $p$, let $\mathbb{Z}_n=\{0,1,\ldots,n-1\}$ and $\mathbb{Z}_{n,p}=\{i\in \mathbb{Z}_n: i\not\equiv 0~({\rm mod}~p)\}$. 
We say that a set $A\subseteq \mathbb{Z}_n$ contains a \emph{cyclic $p$-interval} if there are $a_1,a_2,\ldots,a_p\in A$ such that $a_{i+1}\equiv a_i+1~({\rm mod}~n)$ for $i=1,\ldots,p-1$.
Note that $|\mathbb{Z}_{\Delta+\lceil\frac{\Delta}{\delta-1}\rceil,\delta}|=\Delta$ and
let $X=\{x_i : i\in \mathbb{Z}_\delta\}$ and $Y=\{y_j:j\in \mathbb{Z}_{\Delta+\lceil\frac{\Delta}{\delta-1}\rceil,\delta}\}$.
Then we define an edge colouring $c:E\to \mathbb{Z}_{\Delta+\lceil\frac{\Delta}{\delta-1}\rceil}$ by setting $c(x_iy_j)=i+j~({\rm mod}~ \Delta+\lceil\frac{\Delta}{\delta-1}\rceil)$ for all $x_i\in X$, $y_j\in Y$. 
Note that $c$ is a proper edge colouring of $G$ with (at most) $\Delta+\lceil\frac{\Delta}{\delta-1}\rceil$ colours.
As $\delta\leq \Delta$, in order to prove that $c$ is inclusion-free, it suffices to show that for every $y\in Y$ and each $x\in X$, $S_c(y)\not\subset S_c(x)$. This however follows by definition of $c$ as it implies that $S_c(y)$ contains a cyclic $\delta$-interval $A\subseteq  \mathbb{Z}_{\Delta+\lceil\frac{\Delta}{\delta-1}\rceil}$, while $S_c(x)$ does not.
$\Box$
\end{pf}

\section{Conjecture and Main Result}\label{SectionConjecture}

In some sense bipartite graphs with large degrees in one set of the bipartition and small degrees in the other one seem to be most problematic while designing inclusion-free colourings. This is e.g. the case within the proof of our main result below, where a considerable part of the construction is devoted to analysis of a certain bipartite subgraph of a given graph induced by edges between small and large degree vertices (with a special set of colours dedicated in the main to overcome obstacles concerning adjacencies of such vertices). Taking into account this in combination with Proposition~\ref{ThBipartiteComplIncl} above, we pose the following conjecture, which by the mentioned proposition would thus be sharp if proven. (Note that as our requirement is local, i.e. concerns only adjacent vertices, we may focus on connected graphs.)

\begin{conjecture}\label{MainInclConj} For any connected graph $G$ of minimum degree $\delta\geq 2$ and maximum degree $\Delta$ which is not isomorphic to $C_5$,
$$\chi'_{\subset}(G) \le \left\lceil\left(1+\frac{1}{\delta-1}\right)\Delta\right\rceil.$$
\end{conjecture}

Note that in all cases, $\lceil(1+\frac{1}{\delta-1})\Delta\rceil\geq \Delta+2$, where the equality holds in particular for regular graphs. Since the problems of adjacent vertex distinguishing colourings and inclusion-free colourings are equivalent for such a class of graphs, we immediately obtain that Conjecture~\ref{MainInclConj} holds e.g. for
cubic graphs, cycles, regular bipartite graphs and complete graphs, see~\cite{BalGLS,Zhang}. The main contribution of this paper is the following upper bound for the value of the inclusion chromatic index of a graph.

\begin{theorem} \label{thm:main}
%There exists an absolute constant $C>0$ such that for 
For every fixed $\delta\geq 2$ and a graph $G$ of minimum degree $\delta$ and maximum degree $\Delta$,
$$\chi'_\subset(G)\leq \left(1+\frac{4}{\delta}\right)\Delta + O\left(\Delta^{2/3}\log^4\Delta\right).$$ 
%+ C\Delta^{2/3}\log^4\Delta.$$
\end{theorem}
Note this immediately implies the following corollaries.
\begin{corollary}
For any fixed $\delta\geq 2$ and $\Delta$ large enough, $\chi'_\subset(G)\leq \left(1+\frac{4}{\delta-1}\right)\Delta$ for every graph $G$ with minimum degree $\delta$ and maximum degree $\Delta$.
\end{corollary}
\begin{corollary}
For any fixed $\delta\geq 2$ there exists a constant $C$ such that $\chi'_\subset(G)\leq \left(1+\frac{4}{\delta-1}\right)\Delta+C$ for every graph $G$ with minimum degree $\delta$ and maximum degree $\Delta$.
\end{corollary}

In order to prove Theorem~\ref{thm:main} we shall use probabilistic approach, based in particular on the following variants of the Lov\'asz Local Lemma, see e.g.~\cite{AlonSpencer},
the Chernoff Bound, see e.g.~\cite{JansonLuczakRucinski} (Th. 2.1, page 26)
and Talagrand's Inequality, see e.g.~\cite{MolloyReed_GoodTalagrand}.

\begin{theorem}[\textbf{The Local Lemma}]
\label{LLL-symmetric}
Let $A_1,A_2,\ldots,A_n$ be events in an arbitrary pro\-ba\-bi\-li\-ty space.
Suppose that each event $A_i$ is mutually independent of a set of all the other
events $A_j$ but at most $D$, and that $\mathbf{Pr}(A_i)\leq p$ for all $1\leq i \leq n$. If
$$ ep(D+1) \leq 1,$$
then $ \mathbf{Pr}\left(\bigcap_{i=1}^n\overline{A_i}\right)>0$.
\end{theorem}
\begin{theorem}[\textbf{Chernoff Bound}]\label{ChernofBoundTh}
For any $0\leq t\leq np$,
$$\mathbf{Pr}({\rm BIN}(n,p)>np+t)<e^{-\frac{t^2}{3np}}~~{and}~~\mathbf{Pr}({\rm BIN}(n,p)<np-t)<e^{-\frac{t^2}{2np}}\leq e^{-\frac{t^2}{3np}}$$
where ${\rm BIN}(n,p)$ is the sum of $n$ independent Bernoulli variables, each equal to $1$ with probability $p$ and $0$ otherwise.
\end{theorem}

\begin{theorem}[\textbf{Talagrand's Inequality}]\label{TalagrandsInequalityTotal}
Let $X$ be a non-negative random variable determined by $l$ independent trials $T_1,\ldots,T_l$.
Suppose there exist constants $c,k>0$ such that for every set of possible outcomes of the trials, we have:
\begin{itemize}
\item[1.] changing the outcome of any one trial can affect $X$ by at most $c$, and
\item[2.] for each $s>0$, if $X\geq s$ then there is a set of at most $ks$ trials whose outcomes certify that $X\geq s$.
\end{itemize}
Then for any $t\geq 0$, we have
\begin{equation}\label{TalagrandsInequality}
\mathbf{Pr}(|X-{\mathbf E}(X)|>t+20c\sqrt{k{\mathbf E}(X)}+64c^2k)\leq 4e^{-\frac{t^2}{8c^2k({\mathbf E}(X)+t)}}.
\end{equation}
\end{theorem}

Prior to proving Theorem~\ref{thm:main}, we present a deterministic proof of the fact that $\chi'_{\subset}(G)< 3\Delta$, which shall also be useful in the proof of our main result. Note in this context that the bound from Conjecture~\ref{MainInclConj} (and Proposition~\ref{ThBipartiteComplIncl}) may be as large as $2\Delta$ (for $\delta=2$).

\section{Deterministic Bound}

\begin{theorem} \label{thm:3delta} Let $G=(V,E)$ be a graph of minimum degree $\delta \ge 2$ and maximum degree $\Delta$. Then
$$\chi'_{\subset}(G) \le 3\Delta-1.$$
\end{theorem}
\begin{pf} 
Let $v_1, v_2, \dots, v_n$ be an ordering of the vertex set of $G$ such that each vertex $v_i$ has  minimum degree in $G_i = G[\{v_1, \dots, v_i\}]$. We show that there is a proper edge colouring of $G_i$ with a set $C$ of $3\Delta-1$ colours such that for each $i$, if $xy\in E(G_i)$ and $d_{G_i}(x)\ge2$ then $S_i(x)\not\subset S_i(y)$ where $S_i(x)$ is the palette of $x$ in a colouring of $G_i$. In fact we construct the colouring greedily, and in step $i$ we extend the colouring of $G_{i-1}$ by choosing appropriate colours only for the edges incident with $v_i$ in $G_i$.

For $i=1$ and $i=2$ the claim is trivial; assume thus that $i\ge 3$. 

If $d_{G_i}(v_i)=0$, then we have nothing to do and the claim still holds (as it was satisfied for $G_1,\ldots, G_{i-1}$).

If $d_{G_i}(v_i)=1$ then we choose a colour only for one edge, say $uv_i$. We need to preserve the properness and the distinction between $u$ and its neighbours other than $v_i$. So if there is a colour $\alpha$ in the palette of the neighbour $u'\neq v_i$ of $u$ such that using that colour on the edge $uv_i$ would cause $S_i(u')\subset S_i(u)$, we shall forbid it for $uv_i$. In other words, we define
$$\bar{S}_{i-1}(u) = \left\{\alpha\in C\mid \exists u'\in N_{G_{i-1}}(u) \colon S_{i-1}(u')\setminus S_{i-1}(u) = \{\alpha\}\right\}.$$
If $d_{G_i}(u)\ge3$ then $S_{i-1}(u)\not\subset S_{i-1}(u')$ for any neighbour $u'\neq v_i$ of $u$ by the properties of our colouring of $G_{i-1}$, so $S_i(u)\not\subset S_i(u')$ no matter what we put on $uv_i$. Therefore if we choose any element of $C \setminus (S_{i-1}(u)\cup \bar{S}_{i-1}(u))$ as a colour of $uv_i$, we also guarantee that $S_i(u')\not\subset S_i(u)$ for every neighbour $u'$ of $u$ with $d_{G_i}(u')\geq 2$. On the other hand, if $d_{G_i}(u)=2$ then we choose any element of $C\setminus S_{i-1}(u')$ where $u'$ is the only neighbour of $u$ in $G_{i-1}$. If finally $d_{G_i}(u) = 1$, we may choose any colour in $C$ for $uv_i$.  In each case we succeed as long as $|C|\ge 2\Delta-1$.

Now let us assume that $d_{G_i}(v_i)=k\ge 2$ and let $N_{G_i}(v_i) = \{u_1, \dots, u_k\}$.

By the assumption regarding $v_i$ we have $d_{G_i}(v_i) \le d_{G_i}(u_j)$, $j=1,\dots,k$, so there is only one type of inclusion along this edge we need to deal with. In fact we need to choose pairwise distinct colours $c_1,\dots,c_k\in C$ to be assigned to $v_iu_1,\ldots,v_iu_k$, resp.,  such that for each $j\in \{1,\dots,k\}$:  $c_j\not\in S_{i-1}(u_j) \cup \bar{S}_{i-1}(u_j) $ (or  $c_j\notin S_{i-1}(u_j')$ if $d_{G_i}(u_j)=2$ and $u_j'$ is the only neighbour of $u_j$ in $G_{i-1}$) and $c_r\not\in S_{i-1}(u_j)$ for some $r$ other than $j$. We shall choose such $c_j$'s one after another so that additionally for $j=1,\ldots,k-1$, there is $r\in \{1,\dots,j\}$ such that $c_r\notin S_{i-1}(u_{j+1})$  
(what shall guarantee fulfillment of the last requirement above for $u_{2},\ldots,u_k$).
Suppose we have chosen appropriate $c_j$'s for all $j\leq j'-1$ for some $j'\in \{1,\dots,k\}$, and we are about to choose $c_{j'}$. In the following all indices are regarded modulo $k$.

If all colours $c_1,\dots,c_{j'-1}$ are in $S_{i-1}(u_{j'+1})$ (note we must have $j'<k$ then), then we choose as $c_{j'}$ any element of the set $C\setminus (S_{i-1}(u_{j'}) \cup \bar{S}_{i-1}(u_{j'}) \cup S_{i-1}(u_{{j'}+1}))$ or of the set $C\setminus (S_{i-1}(u'_{j'})  \cup S_{i-1}(u_{{j'}+1}))$ if $d_{G_i}(u_{j'})=2$ and $u'_{j'}$ is the only neighbour of $u_{j'}$ in $G_{i-1}$.

We also proceed in the same way in the case when $j'=k$ and the only $c_j$ outside $S_{i-1}(u_{1})$ is $c_1$ but we additionally cannot use $c_1$ as $c_k$ then.

On the other hand, if there is $c_r$ with $r<j'$ and $r\neq 1$ for $j'=k$ such that $c_r\not\in S_{i-1}(u_{j'+1})$, then we choose as $c_{j'}$ any element of the set $C\setminus (S_{i-1}(u_{j'}) \cup \bar{S}_{i-1}(u_{j'}) \cup \{c_1,\ldots,c_{j'-1}\})$ or of the set $C\setminus (S_{i-1}(u'_{j'})  \cup \{c_1,\ldots,c_{j'-1}\})$ if $d_{G_i}(u_{j'})=2$ and $u'_{j'}$ is the only neighbour of $u_{j'}$ in $G_{i-1}$.

It is straightforward to verify that in every case these choices can be committed if only $|C|\geq 3\Delta-1$, and the resulting colouring fulfills all our requirements.
$\Box$
\end{pf}

\section{Proof of Theorem \ref{thm:main}}\label{ProbabilisticProofInclusion}

For $\delta=2$ the thesis follows by Theorem~\ref{thm:3delta}.

Fix an integer $\delta\ge 3$. Let $G=(V,E)$ be a graph with minimum degree $\delta$ and maximum degree $\Delta$. Some of the claims and explicit inequalities in the following proof are true only for $\Delta$ large enough, which we do not specify, say $\Delta \geq \Delta_0$. 
%where $\Delta_0\geq 4$ is some constant. We shall prove the theorem for $C=\Delta_0$ and any $\Delta\geq %\Delta_0$. Note that for such $C$ the thesis follows by Theorem \ref{thm:3delta} whenever $\Delta<\Delta_0$.
We shall in fact prove that if $\Delta$ is sufficiently large, then 
$\chi'_\subset(G)\leq \Delta(1+\frac4\delta) + 4\Delta^{2/3}\log^4\Delta$, what e.g. due to Theorem~\ref{thm:3delta} implies Theorem~\ref{thm:main}.

In the following, $d(v)$ shall always be understood as the degree of a vertex $v$ in $G$. (i.e., $d(v)=d_G(v)$). We shall colour the edges of $G$ with integer colours $1,2,\ldots$ and we shall denote the colour assigned to a given edge $e$ by $c(e)$.

Let 
$$S=\{v\in V: d(v)\leq \Delta^{2/3}\log^4\Delta\},$$
$$B=\{v\in V: d(v) > \Delta^{2/3}\log^4\Delta\},$$
and let $G_S=G[S]$, $G_B=G[B]$ be the subgraphs induced by $S$ and $B$, resp., in $G$. Set
$$S_0 = \{v\in S: d_{S}(v)=0\},$$
$$S_1 = \{v\in S: d_{S}(v)=1\},$$
and let $H$ be the bipartite graph with the sets of bipartition $X=S_0\cup S_1$ and $Y=B$ induced by all the edges between these two sets in $G$. Note that $3\leq\delta\leq d_H(v)\leq \Delta^{2/3}\log^4\Delta$ for every $v\in S_0$,  $2\leq\delta-1\leq d_H(v)< \Delta^{2/3}\log^4\Delta$ for $v\in S_1$ and  $0\leq d_H(v)\leq \Delta$ for $v\in B$.

\begin{claim} \label{claim-G'} There exists a subgraph $G'$ of $G-E(G_S)$ and its subgraph $H'$ which is also a subgraph of $H$ such that:
\begin{enumerate}[(i)]
\item $d_{H'}(v)= 2$ for $v\in S_0$;
\item $d_{H'}(v)= 1$ for $v\in S_1$;
\item $d_{G'}(v)\geq  \frac2\delta d(v) -2\sqrt\Delta\log\Delta$ for every $v\in B$;	\label{ineq-subg-upper}
\item $d_{G'}(v)\leq d(v)\frac2{\delta}(1+\frac{\delta}{\Delta^{1/3}})$ for every $v\in B$ (hence $\Delta(G')\leq \frac{2\Delta}{\delta}+2\Delta^{2/3}$). \label{ineq-subg-lower}
\end{enumerate} \end{claim}

\begin{pf} For every vertex $v\in S_0$ we independently and equiprobably choose a pair of edges of $H$ incident with $v$ and include them in the randomly constructed $G'$ and $H'$; independently, for every edge incident with such $v$ in $H$ we make an equalizing coin flip and include this edge in $G'$ (if it is not yet included in $G'$), but not in $H'$, with probability $\frac2{\delta}(1-\frac{\delta-2}{d_H(v)-2})$, so that every edge incident with $v\in S_0$ lands in $G'$ with probability $2/\delta$;

Similarly, for every vertex $v\in S_1$ we independently and equiprobably choose one edge of $H$ incident with $v$ and include it in the randomly constructed $G'$ and $H'$; independently, for every edge incident with such $v$ in $H$ we make an equalizing coin flip and include this edge in $G'$ (if it is not yet included in $G'$), but not in $H'$, with probability $\frac2{\delta}(1-\frac{\delta-2}{2d_H(v)-2})$ so that every edge incident with $v\in S_1$ in $H$ (thus now also every edge in $H$) lands in $G'$ with probability $2/\delta$;

Finally, every edge of $G-(E(H)\cup E(G_S))$  is independently included in $G'$ with probability $2/\delta$.

The obtained subgraphs satisfy conditions (i) and (ii) by design. We use Talagrand's inequality to bound the probability for the remaining conditions to hold. Let for every $v\in B$, $A_v$ and $B_v$ denote the events that \eqref{ineq-subg-upper} and \eqref{ineq-subg-lower} do not hold for $v$, respectively. The independent trials here are simply the random choices described above determining whether or not a specific edge is included in $G'$. Note that changing the outcome of any trial can affect $d_{G'}(v)$ by at most 1 and for each $s>0$, if $d_{G'}(v)\geq s$ then there is a set of at most $s$ trials whose outcomes certify that $d_{G'}(v)\geq s$. The expected value of $d_{G'}(v)$ equals $\frac2\delta d(v)$, so by the Talagrand's Inequality, 
\begin{align*} \mathbf{Pr} \left( A_v \right) &\le \mathbf{Pr}\left(d_{G'}(v) < \frac2\delta d(v) - 2\sqrt\Delta\log\Delta\right) 
\le 4e^{-0.5\log^2\Delta} 
%4\exp\left(-\frac1{100} \log^2\Delta\right) 
< \frac1{\Delta^3} 
\end{align*}

Analogously, 
\begin{align*} \mathbf{Pr} &\left( B_v \right) \le \mathbf{Pr} \left( d_{G'}(v) > \frac2{\delta} d(v)  +  \frac{2d(v)}{\Delta^{1/3}}\right) 
\le 4e^{-0.5\log^2\Delta} 
%4\exp\left(-\frac1{100} \log^2\Delta\right) 
< \frac1{\Delta^3} \end{align*}

Therefore, as every event $A_v$ and $B_v$ is mutually independent of all other such events except those associated with vertices $v'$ at distance at most $2$ from $v$, i.e. all except at most $2\Delta^2+1$ other events, the theorem follows by the Lov\'asz Local Lemma.
$\Box$
\end{pf}

\begin{claim}\label{2r-3rColouring} We may colour properly the edges of $G'$ using $2r$ colours from the set $\{1,2,\ldots,3r\}$ with $r= \left\lceil \frac3\delta \Delta(1+\frac{\delta}{\Delta^{1/3}})\right\rceil $ so that: 
\begin{enumerate}[(i)] 
\item if $u\in Y$, then the colours of any two edges incident with $u$ in $G'$ differ by at least $2$; \label{claim-disting-1}
 \item if $v\in S_0$, then some two consecutive integers appear as colours of edges incident with $v$ in $G'$; \label{claim-disting-2}
\item if $v\in S_1$ and $uv\in E(G_S)$ then there is a colour incident with $v$ which is not incident with $u$. \label{claim-disting-3}
\end{enumerate}\end{claim}

\begin{pf} We first modify $G'$ as follows: for every vertex $v\in S_0$ we delete the two edges $vu$ and $vw$ incident in $H'$ with $v$, and replace them with one edge $uw$.  
Note that this way we may obtain a multigraph from $G'$ -- we denote it by $G''$. Observe that by Claim~\ref{claim-G'} (iv), $\Delta(G'')\leq 2\frac\Delta\delta(1+\frac{\delta}{\Delta^{1/3}})$, where by the degree of a vertex $v$ in $G''$ we mean the number of edges incident with $v$ in $G''$. 
Consider pairs of consecutive integers: $P_i:=\{3i-2,3i-1\}$ for $i=1,2,\ldots,\left\lceil 3\frac\Delta\delta(1+\frac{\delta}{\Delta^{1/3}})\right\rceil$. By Shannon's theorem \cite{shannon} we may properly colour the edges of $G''$ with colours $P_1,P_2,\ldots$. We shall use this colouring to obtain a desired proper edge colouring of $G'$. We shall never use colours $3,6,\dots,3r$.

Now we process the edges between $Y$ and $S\setminus S_0$ by looking at each connected component of $G[S\setminus S_0]$ one after another, and choosing colours for all edges incident with it in $G''$. If a given component is an isolated edge $uv$ contained in $S_1$, we choose distinct colours from the pairs $P_i$ assigned to the edges incident with $u$ and $v$. Otherwise we choose smaller colours from all pairs $P_i$ assigned to the edges incident with vertices in $S_1$, and larger for the remaining ones, i.e. those between $Y$ and $S\setminus (S_0\cup S_1)$ from the component.
(As we shall not use the colours from this step in the further part of the construction, this shall guarantee the set of colours incident with any vertex $v\in S_1$ shall not be contained in the set of  colours incident with any its neighbour in $S$.)

Then for every vertex $v\in S_0$, one after another, if the corresponding $uw$ (where $vu,vw\in H'$) has colour $P_i$ and this colour appears on some edge incident with $v$ (in $G''$), say $vy$, then we recolour (greedily) this edge with a different $P_j$ so that the colouring of $G''$ remains proper -- this is always feasible as by Claim~\ref{claim-G'} (iv), $d_{G''}(y)+d_{G''}(v)\leq d_{G'}(y)+d_{G'}(v) \leq 2\frac{\Delta}{\delta}+2\Delta^{2/3}+\Delta^{2/3}\log^4\Delta < 3\frac{\Delta}{\delta}$;  
then we replace back $uw$ with $vu$ and $vw$, and colour them with different elements of $P_i$. Note that in this step we did not change any colours on the edges incident with $S\setminus S_0$. For each remaining edge of $G'$
we can pick any colour from its corresponding $P_i$.
$\Box$
\end{pf}

Note that the conditions \eqref{claim-disting-1} and \eqref{claim-disting-2} in Claim~\ref{2r-3rColouring} guarantee that in the obtained partial edge colouring of $G$ the vertices in $S_0$ are distinguished from their neighbours in $B$. Moreover, the condition~\eqref{claim-disting-3} assures that the vertices in $S_1$ are distinguished from their neighbours in $S$.

We then take two copies of $G_S-S_0$ and add a matching between the corresponding vertices of degree $1$. This graph has minimum degree at least $2$ and maximum degree at most $ \Delta^{2/3}\log^4\Delta$, so by Theorem \ref{thm:3delta} we can properly colour its edges using new colours in $\{3r+1,\ldots,3r+s\}$ for some integer $s\leq 3\Delta^{2/3}\log^4\Delta$ so that back in the original copy of $G_S$ all neighbours in $S$ are distinguished in both directions in $G_S\cup G'$ (by the paragraph above this holds also for vertices of degree $1$ in $G_S$). These distinctions shall not be spoiled as the colours used thus far shall not be repeated in the further part of the construction of edge colouring of the entire $G$.

Note moreover that within this step, every vertex in $S\setminus S_0$ received a colour into its palette that shall never appear in the palette of any vertex in $B$. Therefore the only possible inclusions left are between the vertices of $B$.

Let $G_1$ be the subgraph of $G$ induced by all its yet not coloured edges. Note that by Claim \ref{claim-G'}, for every $v\in B$:
\begin{equation}\label{EqdG1-low} 
d_{G_1}(v) \geq d(v)\left(1-\frac{2}{\delta}-\frac2{\Delta^{1/3}}\right) \geq \frac{d(v)}{4} \geq \frac{\Delta^{2/3}\log^4\Delta}{4} 
\end{equation}
and
\begin{equation}\label{EqdG1-up}
d_{G_1}(v)  \leq d(v)\left(1-\frac{2}{\delta}\right) + 2\sqrt\Delta\log\Delta \le \left(1-\frac{2}{\delta} + \frac{2\log\Delta}{\sqrt\Delta}\right)\Delta.
\end{equation}

\begin{claim} \label{SubgraphF} There is a subgraph $F$ of $G_1$ such that:
\begin{enumerate}[(i)] 
\item $\frac18 \log^2\Delta \le d_F(v) \le 2\frac{\Delta^{1/3}}{\log^2\Delta}$ for each $v\in B$ and\label{SubgraphF1}
\item $d_F(v) \le 2\log^2\Delta$ for each $v\in S$.
\end{enumerate}\end{claim}

\begin{pf} Choose every edge of $G_1$ randomly and independently with probability $\frac1{\Delta^{2/3}\log^2\Delta}$ and denote the subgraph induced by the chosen edges by $F$.  Since for $v\in V(G_1)$, $d_F(v) \sim \BIN(d_{G_1}(v), \frac1{\Delta^{2/3}\log^2\Delta})$, by (\ref{EqdG1-low}) %, (\ref{EqdG1-up}) 
and the Chernoff Bound we have:
\begin{eqnarray} \mathbf{Pr}\left( d_F(v) < \frac18 \log^2\Delta ~~\vee~~ d_F(v) > 2\frac{\Delta^{1/3}}{\log^2\Delta}\right)  &\le& \nonumber\\
\le \mathbf{Pr}\left( d_F(v) < \frac18 \log^2\Delta\right) + \mathbf{Pr}\left( d_F(v) > 2\frac{\Delta^{1/3}}{\log^2\Delta}\right)  &\le& 2e^{-\frac1{32} \log^2\Delta} < \frac1{\Delta^2} \nonumber
 \end{eqnarray}
for $v\in B$ and, by the definition of $S$,
$$ \mathbf{Pr}\left( d_F(v) > 2\log^2\Delta\right)  \le e^{-\frac14 \log^2\Delta} < \frac1{\Delta^2}$$
for $v\in S$.

Any of the analysed events, opposite to the ones from the thesis of Claim~\ref{SubgraphF}, is mutually independent of all the other such events but the ones associated with the vertices at distance at most $1$ from $v$, i.e. all other events but at most $2\Delta+1$. The thesis thus follows by the Lov\'asz Local Lemma.
$\Box$
\end{pf}

Now we randomly and independently colour each edge $e$ of $F$ with a colour $c'(e)\in\{3r+s+1,\ldots,3r+s+t\}$, $t=\lceil\Delta^{1/3}\rceil$ -- choosing every colour with probability $1/\lceil\Delta^{1/3}\rceil\leq\Delta^{-1/3}$. At the end we uncolour every edge $e$ adjacent with an edge $e'$ coloured the same, i.e. with $c'(e)=c'(e')$ (note that $e'$ is uncoloured as well then).

\begin{claim} It is possible to make our choices so that after the uncolouring, for any $uv\in E(G_B)$, there are still at least $\frac14 d_F(u)-1$ new colours in the palette of $u$ that are not present in the palette of $v$.
\end{claim}

\begin{pf} For every $u,v\in B$ such that $uv\in E$ (note we consider neighbours in $G$, not only in $F$), we define the set:
\begin{eqnarray}
I_{u,v}=\{x\in N_F(u)\smallsetminus\{v\}|&& (\forall z\in N_F(u),z\neq x: c'(ux)\neq c'(uz))\nonumber\\ 
&\wedge&
(\forall z\in N_F(x),z\neq u: c'(ux)\neq c'(xz))\nonumber\\ 
&\wedge& 
(\forall y\in N_F(v): c'(vy)\neq c'(ux))\},\nonumber
\end{eqnarray}
whose cardinality simply represents the number of edges incident with $u$ in $F$ which were not uncoloured and whose colours were not assigned to any edge incident with $v$. 
(Note also that $I_{v,u}$ is not the same as $I_{u,v}$.) 
Let $C_{u,v}=(N_F(u)\smallsetminus\{v\})\smallsetminus I_{u,v}$, i.e.,
\begin{eqnarray}
C_{u,v}=\{x\in N_F(u)\smallsetminus\{v\}| && (\exists z\in N_F(u),z\neq x: c'(ux)= c'(uz))\nonumber\\ 
&\vee& (\exists z\in N_F(x),z\neq u: c'(ux)= c'(xz))\nonumber\\ 
&\vee& (\exists y\in N_F(v): c'(vy) = c'(ux))\}.\nonumber
\end{eqnarray}
We shall show the colour choices can be made so that $|I_{u,v}|\geq \frac{1}{4}d_F(u)-1$ for every ordered pair $(u,v)$ such that $u,v\in B$ and $uv\in E$. Thus 
%for such vertices 
let us denote the following opposite event for such $(u,v)$:
$$A_{u,v}: |C_{u,v}| > \frac{3}{4}d_F(u).$$
For any fixed $x\in N_F(u)\smallsetminus\{v\}$, by Claim~\ref{SubgraphF}, %~\eqref{SubgraphF1}, 
we have:
\begin{eqnarray}
\mathbf{Pr}(x\in C_{u,v}) \leq&& \mathbf{Pr} (\exists z\in N_F(u),z\neq x: c'(ux)= c'(uz))\nonumber\\
&+& \mathbf{Pr} (\exists z\in N_F(x),z\neq u: c'(ux)= c'(xz))\nonumber\\
&+& \mathbf{Pr}(\exists y\in N_F(v): c'(vy) = c'(ux))\nonumber\\
\leq&& \sum_{z\in N_F(u),z\neq x} \mathbf{Pr} (c'(ux)= c'(uz)) \nonumber\\
&+& \sum_{z\in N_F(x),z\neq u} \mathbf{Pr} (c'(ux)= c'(xz)) \nonumber\\
&+& \sum_ {y\in N_F(v)} \mathbf{Pr}(c'(vy) = c'(ux))\nonumber\\
\leq&& \left( d_F(u)+d_F(x)+d_F(v)\right) \cdot \Delta^{-1/3} \nonumber\\
\leq&& 3\cdot 2\frac{\Delta^{1/3}}{\log^2\Delta}\cdot \Delta^{-1/3} =\frac{6}{\log^2\Delta}.\nonumber
\end{eqnarray}
Consequently,
$${\mathbf E} \left(|C_{u,v}\right|) \leq d_F(u)\cdot \frac{6}{\log^2\Delta} \nonumber \leq \frac{1}{2}d_F(u).$$
The random variable $|C_{u,v}|$ depends only on the result of the assignment of colours in this step (i.e. analysed within this claim). The change of one colour can change $|C_{u,v}|$ by at most 2. Moreover, for any $s$, we can indicate a set of at most $2s$ edges whose colours certify that $|C_{u,v}|\ge s$. Therefore, by Talagrand's Inequality
and Claim~\ref{SubgraphF}~\eqref{SubgraphF1}:
$$\mathbf{Pr}(A_{u,v}) = \mathbf{Pr}\left(|C_{u,v}|> \frac{3}{4}d_F(u)\right) \leq 4e^{-0.001\cdot d_F(u)} \leq 4e^{-0.0001 \log^2\Delta} \leq \Delta^{-5}.$$
%$$\mathbf{Pr}(A_{u,v}) = \mathbf{Pr}\left(|C_{u,v}|> \frac{3}{4}d_F(u)\right) \leq 4e^{-\frac1{640}d_F(u)} \leq 4e^{-%\frac1{5120} \log^2\Delta} \leq \Delta^{-5}.$$
As every event $A_{u,v}$ (where $A_{u,v}$ is not the same as $A_{v,u}$) is mutually independent of all other such events associated with edges at distance at least $5$ from the edge $uv$, i.e. of all except at most $4\Delta^4$ events, then by the Lov\'asz Local Lemma, with positive probability none of the events $A_{u,v}$ holds. 
$\Box$
\end{pf}

For every not uncoloured edge of $uv\in F$, we set $c(uv)=c'(uv)$. Note that all neighbours in $G$ are distinguished in both directions by the obtained partial colouring with $2r+s+t$ colours.

By~(\ref{EqdG1-up}), the definition of $S$ and Vizing's theorem, the remaining edges may be coloured properly with at most 
$\left(1-\frac{2}{\delta} + \frac{2\log\Delta}{\sqrt\Delta}\right)\Delta+1$
 new colours, what finalizes the construction of a desired proper colouring of $G$ with at most
\begin{align*} 2r+s+t &+ \left(1-\frac{2}{\delta} + \frac{2\log\Delta}{\sqrt\Delta}\right)\Delta+1 \leq \\
&\leq \frac6\delta \left(1+\frac{\delta}{\Delta^{1/3}}\right)\Delta + 2 + 3\Delta^{2/3}\log^4\Delta + \Delta^{1/3}+1+\left(1-\frac{2}{\delta} + \frac{2\log\Delta}{\sqrt\Delta}\right)\Delta +1\\ 
&\leq \Delta\left(1+\frac4\delta\right) + 4\Delta^{2/3}\log^4\Delta \end{align*}
colours. 
$\Box$\\

\noindent \textbf{Acknowledgements:} 
The authors would like to thank G. Simonyi for asking the interesting question, which inspired us to write this paper.\\

\end{document}